# SELECT SETS: RANK AND FILE

By Abba M. Krieger,[1] Moshe Pollak[1] and Ester Samuel-Cahn[2]

*University of Pennsylvania, Hebrew University and Hebrew University*

In many situations, the decision maker observes items in sequence and needs to determine whether or not to retain a particular item immediately after it is observed. Any decision rule creates a set of items that are selected. We consider situations where the available information is the rank of a present observation relative to its predecessors. Certain "natural" selection rules are investigated. Theoretical results are presented pertaining to the evolution of the number of items selected, measures of their quality and the time it would take to amass a group of a given size.

**1. Introduction.** Items (people) are observed sequentially, each receiving a score. The decision whether or not to accept an item must be made on the spot, without possibility of getting back to an item that has been let go. What would constitute a reasonable selection policy?

The issues involved in formulating a policy are the quality and quantity of those selected and the size and rate at which the accepted set grows. Also of import is whether or not the horizon of the pool of items is finite. The scenario we envision is one where observations arrive in random order and their scores are independent and identically distributed, but nothing is otherwise known about their distribution so that, as items are observed, information about the pool of candidates is being gathered. Heuristically, quality of the items selected and speed of selecting items are in conflict with each other. A policy of selecting all items fulfills the need for speed, but the quality will be average. Toward the other extreme, declining to accept any item unless it is better than all of those observed previously will produce

Received March 2005; revised September 2006.

[1]Supported by funds from the Marcy Bogen Chair of Statistics at the Hebrew University of Jerusalem.

[2]Supported by the Israel Science Foundation Grant 467/04.

*AMS 2000 subject classifications.* Primary 62G99; secondary 62F07, 60F15.

*Key words and phrases.* Selection rules, ranks, nonparametrics, sequential observations, asymptotics.

a high-quality set of items, but its rate of growth will be very slow. In this paper we study certain policies that compromise between these two objectives.

The procedure that accepts the first item and subsequently accepts only items which are better than all those observed previously is a well-studied policy (cf. [1] and [6]).

Preater [5] studied a method that prescribes the acceptance of the first observation and subsequently accepts those observations that would improve the average score of those retained. Preater assumed that the scores are exponentially distributed, and derived the asymptotic growth and distribution of the average score after $n$ observations have been retained, as $n \to \infty$. Selection rules with known distribution of the items and inspection cost are considered in [4].

Other problems that have a similar flavor are variations of the secretary problem (cf. [2]), where one samples sequentially from a finite pool until one or several items are retained, after which sampling ceases, with the objective being the maximization of the probability of retaining the best in the pool.

In this paper, we study sequential rules that are based on the ranks of the observations. At every stage, we (re-) rank the observed values from the best to the worst, so that the best has rank 1. We consider procedures that retain the $n$th observation if its rank is low enough relative to the ranks of the previously retained observations.

For the sake of illustration, consider the median rule, which prescribes the retention of the $n$th observation if and only if its rank is lower than the median rank of the observations retained previously (i.e., the median of the retained group will be improved). Regarding the speed at which observations are retained, let $L_n$ be the number of items that are retained after $n$ items are observed. We show that the expected number $E(L_n)$ of observations retained after $n$ have been observed is of order $n^{1/2}$ and that $L_n/n^{1/2}$ converges almost surely to a nondegenerate random variable. Regarding the quality of the retained observations, we show that at least half of the observations retained are the very best of all observed heretofore, and that the average rank of the retained observations and its expectation are of order $n^{1/2} \log n$ (implying that almost all of the retained observations are very good).

The paper is organized as follows. In Section 2 we introduce a general class of rules that are characterized by a criterion that ensures that the probability that item $n+1$ is retained is a simple function of its rank and the number of items $L_n$ that have already been retained. In Section 3 we specialize and consider rules that retain an item if it is among the best $100p$ percent of the items already retained. In Section 4 we show that $E(L_n)/n^p$ converges and that $L_n/n^p$ converges almost surely to a nondegenerate random variable. In Section 5 we find the order of the expected value of the average rank of the observations retained by the $p$-percentile rule and that suitably normalized, the average converges a.s. to a nondegenerate random variable. We end with remarks and conclusions in Section 6.

**2. A class of selection rules.** As stated in the previous section, our focus in this article is on selection rules based on ranks. In this section we introduce a class of selection rules that retain an observation if its rank is "low enough," where the threshold of "low enough" is determined solely by the size of the set of observations already retained. The rationale for this has to do with the trade-off between the quality of retained observations and the speed of their accumulation. Heuristically, the more observations retained, the slower one would go about retaining further observations, so the size of the retained set should be a factor in the selection rule. On the other hand, one's evaluation of the quality of an observation depends on all past observations, not only on those retained so far, and one's expectations regarding future observations is the same irrespective of the quality of those already retained. Therefore, there is good reason to require a selection rule to depend only on the size of the retained set of observations and the rank (among all observations) of the present observation. (As for the desire to "improve" the set of retained observations, heuristically, the quality of the retained set is correlated with its size, so at least qualitatively "improvement" is implicit in retained set size. We examine this more formally in Theorem 2.1 and Remark 2.1.)

Formally, let $X_1, X_2, \ldots$ be a sequence of observations so that any ordering of the first $n$ observations is equally likely. A sufficient condition is that the random variables be exchangeable. A special case that satisfies this assumption is when we have independent and identically distributed (i.i.d.) random variables from a continuous distribution. Let $S_n$ be the set of indices of the retained $X$'s after $n$ items have been observed and let $L_n$ be the size of $S_n$. Let $R_i^n$ be the rank of the $i$th observation from among $X_1, \ldots, X_n$, that is, $R_i^n = \sum_{j=1}^n I\{X_j \leq X_i\}$ where $I\{A\}$ is the indicator function of $A$. Thus, $R_n^n$ is the rank of $X_n$ within the set $\{X_i\}_{i=1}^n$, where without loss of generality we assume that "better" is equivalent to "smaller" so that rank 1 is given to the smallest observation, rank 2 to the second smallest, etc. A selection rule of the type we study is defined by an integer-valued function $r(\cdot)$ on the integers such that the observation $X_n$ will be retained if and only if $R_n^n \leq r(L_{n-1})$. In this article, we assume that the first observation is always kept.

Another feature of a reasonable selection procedure is to require that the function $r(\cdot)$ be locally subdiagonal, that is, $r(a+1) \leq r(a) + 1$. Again, the rationale for this has to do with the trade-off between the quality of retained observations and the speed of their accumulation. [To see this, suppose $a$ observations have been retained after $n$ have been observed. The rank of the next retained observation will not exceed $r(a)$. The rank of the succeeding retained observation will not exceed $r(a+1)$. If $r(a+1) > r(a) + 1$, it would mean that after having retained $a + 1$ observations, one would be willing to settle for an observation of lower quality than the acceptance level after having retained $a$ observations. Although that may be reasonable in a case that a quota has to be filled and the pool of applicants is finite, that is not the case we regard here.]

We summarize the above in the following definition.

DEFINITION 2.1. A locally subdiagonal rank selection scheme ($LsD$) is a rule determined by an integer-valued function $r(\cdot)$ with the following properties:

  (i) $r$ is nondecreasing.
  (ii) $r(0) = 1$ and $L_0 = 0$.
  (iii) $r$ is locally subdiagonal, that is, $r(a + 1) \leq r(a) + 1$.
  (iv) For $n \geq 1$, $X_n$ is retained if and only if $R_n^n \leq r(L_{n-1})$. (This implies that the first observation is retained.)

This class contains many rules that make heuristic sense. For instance, the median rule is a $LsD$ rule with $r(1) = r(2) = 1$, $r(3) = r(4) = 2$, and generally $r(2j - 1) = r(2j) = j$. A class of $LsD$ rules is "$k$-record rules." For a fixed value $k$, let $r(j) = \min\{j + 1, k\}$. For $k = 1$, this is the classical record rule, where an element is retained if and only if it is better than all previous observations. "$k$-record rules" have been studied extensively (cf. [3] or [6] and many subsequent papers).

The following theorem is trivial for "$k$-record rules" but is true for any $LsD$ rule. It attests to the high quality of the set of observations retainable by a $LsD$ rule.

Let $N \geq 1$ be any integer either predetermined or random. For example, $N$ can be a stopping rule. A special case of interest is inverse sampling, where the objective is to collect a group of some fixed size $m$, so that

$$N = \inf\{n : L_n = m\}.$$

THEOREM 2.1. *Consider a $LsD$ rule defined by $r(\cdot)$. The $r(L_N)$ best observations among $X_1, X_2, \ldots, X_N$ belong to $S_N$.*

PROOF. Let $X_m$ be the $t$th best observation among $X_1, X_2, \ldots, X_N$ with $t \leq r(L_N)$. Let $a$ be the number of observations among $X_1, X_2, \ldots, X_{m-1}$ that are better than $X_m$ and let $b$ be the number of observations among $X_{m+1}, \ldots, X_N$ that are better than $X_m$. Clearly, $a + b = t - 1$.

If $m \notin S_N$, then the next items retained after iteration $m$ must be better than $X_m$. This implies that $L_N \leq L_{m-1} + b$. Hence

$$r(L_N) \leq r(L_{m-1} + b) \leq r(L_{m-1}) + b.$$

But $m \notin S_N$ implies that $a + 1 > r(L_{m-1})$. Since by assumption $t \le r(L_N)$,

$$a + b + 1 = t \le r(L_N) \le r(L_{m-1}) + b$$

so that $a + 1 \le r(L_{m-1})$, which contradicts the inequality two lines above. $\square$

REMARK 2.1. Because of Theorem 2.1, implicit in the definition of a $LsD$ rule is that it "improves" the retained set. For example, when applying the median rule, the median of the retained set gets better, something that is not transparent when regarding the median rule via its $LsD$ definition. Theorem 2.1 is a more formal presentation of the heuristic stated in the beginning of this section, that the quality of the retained set is correlated with its size, and a $LsD$ rule embodies all three heuristics: (i) the larger the retained set, the slower one goes about retaining more observations; (ii) perception of quality is founded on all previous observations; (iii) one only retains items that "improve" the retained set. Theorem 2.1 means that there is no contradiction between the third heuristic and selection based on the size of the retained set only.

A natural representation of the quality of the group of observations kept (when retention is by ranks) is the average rank of the observations retained. Denote by $Q_n$ the sum of the ranks of the retained set after $n$ observations have been made, so that the average rank $A_n$ is $Q_n/L_n$.

LEMMA 2.1. Let $\mathcal{F}_n$ be the $\sigma$-field generated by the ranks of the i.i.d. continuous random variables $X_1, \ldots, X_n$. Let $Q_n = \sum_{i \in S_n} R_i^n$. The conditional expected behavior of the quantities $L_{n+1}, Q_{n+1}$ and the average rank $A_{n+1}$ given the past, in terms of the corresponding quantities for $n$, results in:

(i) $E(L_{n+1}|\mathcal{F}_n) = L_n + \frac{r(L_n)}{n+1}$,

(ii) $E(Q_{n+1}|\mathcal{F}_n) = \frac{n+2}{n+1}Q_n + \frac{r(L_n)(r(L_n)+1)}{2(n+1)}$,

(iii) $E(A_{n+1}|\mathcal{F}_n) = A_n(1 + \frac{1+L_n-r(L_n)}{(n+1)(L_n+1)}) + \frac{r(L_n)(r(L_n)-1)}{2(n+1)L_n}$.

PROOF. (i) It follows that $L_{n+1}|L_n = L_n + B(\frac{r(L_n)}{n+1})$ where $B(x)$ is a Bernoulli random variable with probability $x$. Hence (i) follows by taking conditional expectations on both sides.

(ii) The sequence $\{Q_n\}$ is nondecreasing. Its growth can be described as follows. If $X_{n+1}$ is retained, and it has rank $k$ among the items retained, then $X_{n+1}$ adds $k$ to the sum of the ranks and 1 for each observation that is inferior to it. Hence, $Q_{n+1} = Q_n + k + [L_n - (k-1)] = Q_n + L_n + 1$. When $X_{n+1}$ is not retained, then the rank of some of the lower-quality retained observations can increase (by 1, if $X_{n+1}$ has lower rank). Note that the distribution of the rank of $X_{n+1}$ (conditional on $\mathcal{F}_n$ and its not being retained) is uniform over $r(L_n) + 1, \ldots, n+1$. Therefore, for $n \geq 1$,

$$
\begin{aligned}
E(Q_{n+1}|\mathcal{F}_n) \\
&= (Q_n + L_n + 1)\frac{r(L_n)}{n+1} \\
&\quad + \frac{n+1-r(L_n)}{n+1}\left[Q_n + E\left(\sum_{\{i \in S_n\}} I\{R_{n+1}^{n+1} < R_i^{n+1}\}|\mathcal{F}_n, R_{n+1}^{n+1} > r(L_n)\right)\right] \\
&= Q_n + (L_n + 1)\frac{r(L_n)}{n+1} \\
&\quad + \frac{Q_n - r(L_n)(r(L_n)+1)/2 - r(L_n)(L_n - r(L_n))}{n+1} \\
&= \frac{n+2}{n+1}Q_n + \frac{r(L_n)(r(L_n)+1)}{2(n+1)}.
\end{aligned}
$$

(iii) If $X_{n+1}$ is retained, then $L_{n+1} = L_n + 1$ and if $X_{n+1}$ is not retained, $L_{n+1} = L_n$. Therefore, the same argument as in the proof of (ii) leads to

$$
\begin{aligned}
E(A_{n+1}|\mathcal{F}_n) \\
&= \frac{Q_n + L_n + 1}{L_n + 1} \cdot \frac{r(L_n)}{n+1} \\
&\quad + \frac{n+1-r(L_n)}{n+1} \\
&\qquad \times \left[\frac{Q_n + E(\sum_{\{i \in S_n\}} I\{R_{n+1}^{n+1} < R_i^{n+1}\}|\mathcal{F}_n, R_{n+1}^{n+1} > r(L_n))}{L_n}\right] \\
&= A_n\left[\frac{L_n}{L_n+1} \cdot \frac{r(L_n)}{n+1}\right] + \frac{r(L_n)}{n+1} \\
&\quad + \frac{n+1-r(L_n)}{(n+1)L_n}\left[Q_n + E\left[\frac{\sum_{\{i \in S_n, R_{n+1}^{n+1} < R_i^{n+1}\}}(R_i^{n+1} - r(L_n))}{n+1-r(L_n)}\Big|\mathcal{F}_n\right]\right] \\
&= A_n\left[\frac{r(L_n)L_n}{(n+1)(L_n+1)} + \frac{n+1-r(L_n)}{n+1}\right] \\
&\quad + \frac{r(L_n)}{n+1} + \frac{1}{(n+1)L_n}\left[Q_n - \frac{r(L_n)(r(L_n)+1)}{2} - r(L_n)(L_n - r(L_n))\right] \\
&= A_n\left[1 + \frac{1+L_n-r(L_n)}{(n+1)(L_n+1)}\right] + \frac{r(L_n)(r(L_n)-1)/2}{(n+1)L_n}.
\end{aligned}
$$

$\square$

Since Lemma 2.1 is central to the derivations in the sequel, we will henceforth assume without loss of generality that the random variables $X_1, X_2, \ldots$ are i.i.d. and continuous.

**3. Percentile rules.** In the following sections, we consider rules that retain items if the item is among the best $100p$ percent among those items that have already been retained.

DEFINITION 3.1. A $p$-percentile rule, for $p$ fixed $(0 < p \le 1)$ is a $LsD$ rule with $r(k) = \lceil pk \rceil$ for $k \ge 1$, where $\lceil x \rceil$ is the smallest integer that is greater than or equal to $x$. Thus, the $n$th item is retained if and only if its rank satisfies $R_n^n \le \lceil pL_{n-1} \rceil$.

To see that the $p$-percentile rule is a $LsD$ rule, note that $\lceil p(a+1) \rceil = \lceil pa + p \rceil \le \lceil pa + 1 \rceil = \lceil pa \rceil + 1$.

REMARK 3.1. Note that the $p$-percentile rule is meaningful even when $p = 1$. In that case, the first observation is kept. The second is kept if it is better than the first observation. In general, an item is kept if it is better than the worst item that has already been retained. It is easy to see that when $p = 1$, $E(L_n|L_{n-1}) = L_{n-1} + L_{n-1}/n$. It is straightforward to show that $E(L_n) = \frac{n+1}{2}$. Hence $E(L_n)/n \to 1/2$. Also, since $E(L_n|L_{n-1}) = \frac{n+1}{n}L_{n-1}$, it follows that $L_n/(n+1)$ is a bounded positive martingale, and therefore converges almost surely. Since the worst item that has already been retained is obviously $X_1$, it follows that $L_n/n$ is asymptotically $U(0,1)$.

**4. Results for the number of retained items.** In this section, we study the behavior of the number of items that are retained after $n$ items are observed, $L_n$, for $p$-percentile rules. It turns out that $L_n$ is of order $n^p$. Hence we consider the quantity $L_n/n^p$. We first show that the expectation of this quantity converges to a finite limit. We then show that this quantity itself converges almost surely to a nondegenerate random variable.

The first result we present is that $E(L_n)/n^p \to c_p$ as $n \to \infty$. For example, this result says that the rule that retains items if they are superior to the median of all items already retained, will be keeping on the order of $\sqrt{n}$ items on the average. The constant $c_p$ depends on

$$d_n \equiv E(\lceil pL_n \rceil - pL_n).$$

The relationship between $c_p$ and $d_1, d_2, \ldots$ is complicated because it depends on all of the $d_j$. It seems impossible to determine $c_p$ analytically, except for $p = 1$, as done in Remark 3.1.

The result, however, only requires that we show that $d_n$ is bounded away from zero. This result is intuitive. For the median rule ($p = 1/2$), $d_n$ is simply $P(L_n \text{ is odd})/2$. Logically, we would expect (it turns out to be justified by empirical analysis) that $P(L_n \text{ is odd}) \to 1/2$ as $n \to \infty$. This is not easy to prove. Similarly, if $p = 1/4$, then $\lceil pL_n \rceil - pL_n$ is either 0, 3/4, 1/2 or 1/4 depending on whether $L_n(\text{mod } 4)$ is $j$ for $j = 0, 1, 2$ or 3, respectively. Since logically each of the four cases should be equally likely (again this appears to be the case by computer analysis), we would anticipate that $d_n \to 3/8$. [We conjecture that if $p$ is an irrational number, then $\lceil L_n p \rceil - L_n p$ converges to $U(0, 1)$ which implies that $d_n \to 1/2$.]

The following lemma shows that $d_n$ is bounded away from zero.

LEMMA 4.1.  *Let $0 < p < 1$  be fixed, and  $\varepsilon = \varepsilon_p = \min\{\frac{p}{2}, \frac{1-p}{2}\}$.  Then $d_n \geq \varepsilon/3$ for all $n$.*

PROOF.  Let $S_\varepsilon = \{j \mid \lceil pj \rceil - pj \leq \varepsilon\}$. Note that if $j \in S_\varepsilon$, then

- $j - 1 \notin S_\varepsilon$. This follows since $\varepsilon + p < 1$, thus $\lceil p(j-1) \rceil = \lceil pj \rceil$. But then $\lceil p(j-1) \rceil - p(j-1) = \lceil pj \rceil - pj + p \geq p > \varepsilon$.
- $j + 1 \notin S_\varepsilon$. This follows since $p - \varepsilon > 0$, thus $\lceil p(j+1) \rceil = \lceil pj \rceil + 1$. Hence $\lceil p(j+1) \rceil - p(j+1) = \lceil pj \rceil - pj + 1 - p > \varepsilon$.

We will show that for all $n \geq 2$ and all $j = 1, 2, \ldots,$

$$(4.1) \qquad P(L_n = j + 1) + P(L_n = j - 1) - P(L_n = j) \geq 0.$$

This will yield the lemma since clearly (4.1) implies $\sum_{j \in S_\varepsilon} P(L_n = j) \leq 2 \times \sum_{j \notin S_\varepsilon} P(L_n = j)$, which in turn implies that $\sum_{j \notin S_\varepsilon} P(L_n = j) \geq 1/3$ so $d_n \geq \varepsilon/3$. Note that (4.1) is trivial for $j > n$.

We prove (4.1) by induction. For $n = 2$ and all $0 < p < 1$ we have $P(L_2 = 1) = P(L_2 = 2) = 1/2$. Thus (4.1) holds for $j = 1, 2$ and $n = 2$.

Now assume (4.1) holds for $2, 3, \ldots, n - 1$. We shall show it holds for $n$. Consider first the values of $j$ for which

$$(4.2) \qquad 2\lceil pj \rceil / n \leq 1.$$

Clearly

$$(4.3) \qquad P(L_n = j - 1) \geq (1 - \lceil p(j-1) \rceil / n) P(L_{n-1} = j - 1),$$

$$(4.4) \qquad \begin{aligned} P(L_n = j) &= (\lceil p(j-1) \rceil / n) P(L_{n-1} = j - 1) \\ &\quad + (1 - \lceil pj \rceil / n) P(L_{n-1} = j), \end{aligned}$$

$$(4.5) \qquad \begin{aligned} P(L_n = j + 1) &= (\lceil pj \rceil / n) P(L_{n-1} = j) \\ &\quad + (1 - \lceil p(j+1) \rceil / n) P(L_{n-1} = j + 1). \end{aligned}$$

Thus

$$P(L_n = j + 1) + P(L_n = j - 1) - P(L_n = j)$$
$$\geq P(L_{n-1} = j + 1) + P(L_{n-1} = j - 1)$$
$$- P(L_{n-1} = j) + 2(\lceil pj \rceil / n) P(L_{n-1} = j)$$
$$- 2(\lceil p(j-1) \rceil / n) P(L_{n-1} = j - 1)$$
$$- (\lceil p(j+1) \rceil / n) P(L_{n-1} = j + 1).$$

(4.6)

However, $\lceil p(j-1) \rceil \leq \lceil pj \rceil$ and $\lceil p(j+1) \rceil \leq 2\lceil pj \rceil$ as $\lceil pj \rceil \geq 1$. Hence, the right-hand side of (4.6) is greater than or equal to

$$(1 - 2\lceil pj \rceil / n)[P(L_{n-1} = j + 1) + P(L_{n-1} = j - 1) - P(L_{n-1} = j)]$$
$$\geq 0$$

(4.7)

where the last inequality in (4.7) follows from (4.2) and the induction hypothesis.

Now consider values of $j$ (if such exist) for which

$$2\lceil pj \rceil / n > 1.$$

(4.8)

Then clearly $j > 1$. Replace (4.3) by

$$P(L_n = j - 1) = (1 - \lceil p(j-1) \rceil / n) P(L_{n-1} = j - 1)$$
$$+ (\lceil p(j-2) \rceil / n) P(L_{n-1} = j - 2)$$

(4.9)

and replace (4.5) by

$$P(L_n = j + 1) \geq (\lceil pj \rceil / n) P(L_{n-1} = j).$$

(4.10)

Then by (4.9), (4.4) and (4.10), it follows that

$$P(L_n = j + 1) + P(L_n = j - 1) - P(L_n = j)$$
$$\geq (2\lceil pj \rceil / n - 1) P(L_{n-1} = j)$$
$$+ (\lceil p(j-2) \rceil / n) P(L_{n-1} = j - 2)$$
$$- (2\lceil p(j-1) \rceil / n - 1) P(L_{n-1} = j - 1).$$

(4.11)

If $2\lceil p(j-1) \rceil / n \leq 1$, then by (4.8) clearly the value in the right-hand side of (4.11) is nonnegative. If

$$2\lceil p(j-1) \rceil / n - 1 > 0,$$

(4.12)

we shall show that (4.12) implies

$$\lceil p(j-2) \rceil / n \geq 2\lceil p(j-1) \rceil / n - 1$$

(4.13)

so that the right-hand side of (4.11) is greater than or equal to

$$2\lceil p(j-1)\rceil/n - 1)[P(L_{n-1}=j) + P(L_{n-1}=j-2) - P(L_{n-1}=j-1)] \geq 0$$

where the last inequality follows from (4.12) and the induction hypothesis. To see (4.13) note that $\lceil p(j-2)\rceil \geq \lceil p(j-1)\rceil - 1$. Thus (4.13) will follow if we show that $(\lceil p(j-1)\rceil - 1)/n \geq 2\lceil p(j-1)\rceil/n - 1$ which is equivalent to

$$(4.14) \qquad\qquad n-1 \geq \lceil p(j-1)\rceil.$$

Since $j \leq n$ are the only values of interest, we have $j-1 \leq n-1$, for which (4.14) clearly holds.  □

We now turn to the main result of showing that the average number of items that are retained is of order $n^p$. From Lemma 2.1(i),

$$E(L_n|L_{n-1}) = L_{n-1} + \lceil pL_{n-1}\rceil/n.$$

Hence,

$$E(L_n|L_{n-1}) = L_{n-1} + pL_{n-1}/n + (\lceil pL_{n-1}\rceil - pL_{n-1})/n.$$

Let $M_n = E(L_n)$. Then

$$(4.15) \qquad\qquad M_n = M_{n-1}(1 + p/n) + d_{n-1}/n.$$

We are now prepared to state and prove:

THEOREM 4.1.  *Let* $0 < p \leq 1$. $E(L_n)/n^p \to c_p$ *as* $n \to \infty$ *with* $0 < c_p < \infty$.

PROOF.  By Remark 3.1, $c_1 = 1/2$. Thus consider a fixed $p$, $0 < p < 1$, and let $T_n = M_n/n^p$. From (4.15) we have that

$$(4.16) \qquad T_n = ((n-1)/n)^p(1 + p/n)T_{n-1} + d_{n-1}/n^{1+p}.$$

The key to the proof is showing that $\Delta_n \equiv T_n - T_{n-1}$ eventually becomes positive and remains positive. Since $T_n = \sum_{j=1}^n \Delta_j$ with $T_0 \equiv 0$ and $T_n$ will be shown to be bounded, it follows that $T_n$ converges.

By the definition of $\Delta_j$ and (4.16),

$$(4.17) \qquad\qquad \Delta_j = b_j T_{j-1} + d_{j-1}/j^{1+p}$$

where $b_j = ((j-1)/j)^p(1 + p/j) - 1$.

The basis of the proof is in the result that $b_j < 0$ and increases to 0 as $j \to \infty$. This is a straightforward calculus argument.

Let $x = 1/j$ and $f(x) = (1-x)^p(1+px) - 1$. Thus, $f(0) = 0$. Also, $f'(x) < 0$ by routine calculus.

From $b_j < 0$ and (4.17) it follows that

$$(4.18) \qquad T_n \leq 1 + \sum_{j=2}^{n} \frac{1}{j^{1+p}} \leq 1 + \int_{x=1}^{n} (1/x)^{1+p} \, dx < 2/p$$

so $T_n$ is bounded. To show that $\Delta_n$ is eventually nonnegative note that [by (4.18)] $\Delta_n \geq 0 \leftrightarrow T_{n-1} \leq -\frac{d_{n-1}}{n^{1+p}b_n}$. It is again a straightforward calculus argument to show that $-\frac{1}{n^{1+p}b_n} \to \infty$. Since by (4.18) $T_n < 2/p$, for all $n$, that coupled with Lemma 4.1 will complete the proof. Consider

$$-j^{1+p}b_j = \left[ 1 - \left( \frac{j-1}{j} \right)^p (1 + p/j) \right] j^{1+p}.$$

Again, let $x = 1/j$ and so $-j^{1+p}b_j$ becomes

$$g(x) = [1 - (1-x)^p(1+px)]/x^{1+p}.$$

We need to show that $g(x) \to 0$ as $x \to 0$. This follows easily by l'Hospital's rule. $\square$

We just showed that $E(L_n/n^p)$ converges as $n \to \infty$. Next we show that $L_n/n^p$ has an almost sure limit. We prove this by showing that $L_n/(n+1)^p$ is a (positive) submartingale and that $E(L_n^2/n^{2p})$ is bounded.

THEOREM 4.2. $\lim_{n \to \infty} E(L_n^2/n^{2p})$ exists and is finite.

PROOF. Let $U_n = E(L_n^2/n^{2p})$. We first show that there exist constants $0 < c_1(p) < c_2(p) < \infty$ such that for all $n \geq 1$

$$(4.19) \qquad c_1(p) < U_n < c_2(p).$$

The left-hand side of inequality (4.19) follows trivially from Theorem 4.1, since $U_n \geq (E(L_n)/n^p)^2 \to c_p^2$. For the right-hand side inequality of (4.19), note that

$$
\begin{aligned}
E(L_n^2 | \mathcal{F}_{n-1}) &= L_{n-1}^2 \left( 1 - \frac{\lceil pL_{n-1} \rceil}{n} \right) + (L_{n-1}+1)^2 \frac{\lceil pL_{n-1} \rceil}{n} \\
(4.20) \qquad &\leq L_{n-1}^2 + (2L_{n-1}+1)\frac{pL_{n-1}+1}{n} \\
&= L_{n-1}^2 \left( 1 + \frac{2p}{n} \right) + L_{n-1}\frac{p+2}{n} + \frac{1}{n}.
\end{aligned}
$$

Thus

$$U_n \leq U_{n-1} \left( \frac{n-1}{n} \right)^{2p} \left( 1 + \frac{2p}{n} \right) + \frac{E(L_{n-1})(p+2)}{n^{1+2p}} + \frac{1}{n^{1+2p}}.$$

Therefore,

$$(4.21) \quad U_n - U_{n-1} < U_{n-1} \left\{ \left(\frac{n-1}{n}\right)^{2p} \left(1 + \frac{2p}{n}\right) - 1 \right\}$$

$$+ \frac{E(L_{n-1}/(n-1)^p)(p+2)}{n^{1+p}} + \frac{1}{n^{1+2p}}.$$

Note that since $f(x) = (1-x)^{2p}(1+2px) - 1$ satisfies $f(0) = 0$ and $f'(x) < 0$ for $x > 0$, it follows that $\left(\frac{n-1}{n}\right)^{2p} \left(1 + \frac{2p}{n}\right) - 1 < 0$. Since $E(L_{n-1}/(n-1)^p)$ is bounded, it follows from (4.21) that (with $U_0 = 0$)

$$U_n = \sum_{j=1}^{n} (U_j - U_{j-1}) < \sum_{j=1}^{\infty} \frac{const}{j^{1+p}} < \infty,$$

which accounts for (4.19).

Now denote $\Delta_j = U_j - U_{j-1}$, so that $U_n = \sum_{j=1}^{n} \Delta_j$. By virtue of (4.19) to complete the proof it suffices to show that $\Delta_j > 0$ for all $j$ sufficiently large. By (4.20),

$$E(L_n^2|\mathcal{F}_{n-1}) = L_{n-1}^2 + 2(L_{n-1} + 1)\frac{\lceil pL_{n-1}\rceil}{n}$$

$$\geq L_{n-1}^2 + (2L_{n-1} + 1)\frac{pL_{n-1}}{n}.$$

Thus

$$\Delta_j \geq U_{j-1} \left\{ \left(\frac{j-1}{j}\right)^{2p} \left(1 + \frac{2p}{j}\right) - 1 \right\} + \frac{pE(L_{j-1}/j^p)}{j^{1+p}}.$$

Now for some $0 < \theta < 1$, by Taylor's theorem,

$$\left(\frac{j-1}{j}\right)^{2p} = \left(1 - \frac{1}{j}\right)^{2p} = 1 - \frac{2p}{j} + \frac{p(2p-1)}{j^2}\left(1 - \frac{\theta}{j}\right)^{-2(1-p)}.$$

Hence there exists a constant $c_3(p) > 0$ such that for all $j \geq 1$

$$\left\{ \left(\frac{j-1}{j}\right)^{2p} \left(1 + \frac{2p}{j}\right) - 1 \right\} > -\frac{c_3(p)}{j^2}.$$

Also, there exists a constant $c_4(p) > 0$ such that $E(L_{j-1}/j^p) > c_4(p)$ for all $j > 1$. But then, for all $j$ sufficiently large, $\Delta_j \geq \frac{-c_3(p)c_2(p) + pc_4(p)j^{1-p}}{j^2} > 0$. $\square$

COROLLARY 4.1.  $\lim_{n\to\infty} \mathrm{Var}(L_n/n^p)$ exists and is finite.

Let $j_n \equiv \lceil pL_n \rceil$ be the cutoff rank such that the $(n+1)$st item is retained if and only if its rank from among the first $n + 1$ observations is less than or equal to $j_n$.

THEOREM 4.3. $L_n/(n+1)^p$ *is a submartingale that converges almost surely as $n \to \infty$ to a nondegenerate finite random variable $\Lambda$ such that $\lim_{n\to\infty} E(L_n/(n+1)^p) = E\Lambda = c_p$, for all $0 < p \le 1$.*

PROOF. Since $j_n = \lceil pL_n \rceil$,

$$E(L_n|\mathcal{F}_{n-1}) = \frac{j_{n-1}}{n}(L_{n-1}+1) + \left(1 - \frac{j_{n-1}}{n}\right)L_{n-1}$$

$$= L_{n-1}\left(1 + \frac{p}{n}\right) + \frac{j_{n-1} - pL_{n-1}}{n}$$

so

$$E\left(\frac{L_n}{(n+1)^p} \,\Big|\, \mathcal{F}_{n-1}\right) = \frac{L_{n-1}}{n^p}\left(\frac{n}{n+1}\right)^p\left(1 + \frac{p}{n}\right) + \frac{j_{n-1} - pL_{n-1}}{n(n+1)^p}$$

$$\ge \frac{L_{n-1}}{n^p}\left[\left(\frac{n}{n+1}\right)^p\left(1 + \frac{p}{n}\right)\right]$$

$$\ge \frac{L_{n-1}}{n^p}.$$

Therefore, $L_n/(n+1)^p$ is a positive submartingale. Because $E(L_n/(n+1)^p)$ and $E(L_n^2/(n+1)^{2p})$ are both bounded (by virtue of Theorems 4.1 and 4.2), Theorem 4.3 follows from the submartingale convergence theorem. $\square$

**5. The quality of the retained group of observations acquired by a $p$-percentile rule.** In the previous sections, the focus was on the size of the group retained by the $p$-percentile rule. Here, attention is focused on its quality.

In general $p$-percentile rules yield a qualitative crop. A prime indication of this is Theorem 2.1—after $n$ observations of which $L_n$ have been retained, the best $\lceil pL_n \rceil$ of all $n$ observations seen heretofore are among the retained set. As will be shown below in this section, the other retained observations are generally also of high quality.

To this end, the following theorem considers the average rank of the retained items $A_n$, which equals $Q_n/L_n$.

THEOREM 5.1. *There exist constants $0 < b_p < \infty$ such that for $0 < p \le 1$, $E(A_n)/a_n(p) \xrightarrow[n\to\infty]{} b_p$, where*

$$a_n(p) = \begin{cases} n^{1-p}, & \text{if } p < 1/2, \\ n^{1/2}\log n, & \text{if } p = 1/2, \\ n^p, & \text{if } p > 1/2, \end{cases}$$

*and*

$$b_p = \begin{cases} c_{1/2}/8, & \text{if } p = 1/2, \\ \dfrac{p^2}{2(2p-1)} c_p, & \text{if } p > 1/2, \end{cases} \quad \text{where } c_p \text{ is the limit of } E(L_n/n^p).$$

PROOF.  From Lemma 2.1(iii), with $r(L_n) = j_n = \lceil pL_n \rceil$,

$$(5.1) \qquad E(A_{n+1}|\mathcal{F}_n) = A_n\left[1 + \frac{1 + L_n - j_n}{(n+1)(L_n + 1)}\right] + \frac{j_n(j_n - 1)/2}{(n+1)L_n}.$$

Let $Y_n = A_n/n^{1-p}$. Equation (5.1) implies

$$(5.2) \qquad E(Y_{n+1}|\mathcal{F}_n) = G_n Y_n + B_n$$

where

$$G_n = \left(\frac{n}{n+1}\right)^{1-p}\left[1 + \frac{1 + L_n - j_n}{(n+1)(L_n + 1)}\right]$$

and

$$(5.3) \qquad B_n = \frac{j_n(j_n - 1)/2}{L_n(n+1)^{2-p}}.$$

We consider $B_n$ first. Since $pL_n \le j_n < pL_n + 1$,

$$\frac{(p^2 L_n - p)/2}{(n+1)^{2-p}} \le B_n < \frac{(p^2 L_n + p)/2}{(n+1)^{2-p}}.$$

By Theorem 4.1, $E(L_n/n^p) \underset{n\to\infty}{\longrightarrow} c_p$, which implies

$$(5.4) \qquad E(B_n)n^{2-2p} \underset{n\to\infty}{\longrightarrow} p^2 c_p/2.$$

We consider $G_n$ next. Let $e_n = pL_n + p - \lceil pL_n \rceil$, so that

$$\frac{1 + L_n - j_n}{L_n + 1} = 1 - p + \frac{e_n}{L_n + 1}.$$

Since $(\frac{n}{n+1})^{1-p} = 1 - \frac{1-p}{n+1} + O(\frac{1}{n^2})$ and since $|e_n| \le 1$,

$$(5.5) \qquad G_n = 1 + \frac{e_n}{(L_n + 1)(n + 1)} + O\left(\frac{1}{n^2}\right).$$

Substituting (5.5) into (5.2) yields

$$(5.6) \qquad E(Y_{n+1}|\mathcal{F}_n) = Y_n + \left[\frac{e_n}{(L_n + 1)(n + 1)} + O\left(\frac{1}{n^2}\right)\right]Y_n + B_n.$$

After taking expectations in (5.6), it follows that

$$(5.7) \quad E(Y_{n+1}) = \sum_{m=0}^{n}\left[E(Y_{m+1}) - E(Y_m)\right] = \sum_{m=1}^{n} E(D_m) + \sum_{m=1}^{n} E(B_m)$$

where $D_m = [\frac{e_m}{(L_m+1)(m+1)} + O(\frac{1}{m^2})]Y_m$ and $Y_0 = 0$.

Our aim is to show that $\sum_{m=1}^{n} E(D_m)$ and $\sum_{m=1}^{n} E(B_m)$ (or variants thereof for $p \geq \frac{1}{2}$) have finite limits as $n \to \infty$. For the first sum, since $L_m \leq m$, it is sufficient to show that $E(\sum_{m=1}^{n} \frac{Y_m}{(L_m+1)(m+1)})$ has a finite limit. Now

$$E\left(\frac{Y_m}{(L_m+1)(m+1)}\right) = E\left(\frac{A_m}{m^{1-p}(m+1)(L_m+1)}\right)$$

$$\leq \frac{1}{m^{2-p}}\left\{E\left(A_m \frac{1}{m^\varepsilon}I\{L_m \geq m^\varepsilon\}\right) + E(A_m I\{L_m < m^\varepsilon\})\right\}.$$

By virtue of Lemma A.1 (in the Appendix) there exists a constant $0 < c_{\varepsilon,p} < \infty$ such that, for $0 < \varepsilon < 1/2$,

$$P(L_m < m^\varepsilon) \leq \frac{c_{\varepsilon,p}}{m^{1-\lceil 1+1/p \rceil \varepsilon}} \qquad \text{for all } 1 \leq m < \infty.$$

Note that $A_m \leq m$. Therefore choosing $0 < \varepsilon < (1-p)/\gamma_p$ (with $\gamma_p = \lceil 1 + \frac{1}{p} \rceil$), it follows that

$$E\left|\frac{e_m Y_m}{(L_m+1)(m+1)}\right| < \frac{1}{m^{2-p+\varepsilon}}E(A_m) + c_{\varepsilon,p}/m^{2-p-\gamma_p \varepsilon}.$$

We now divide the proof into three cases.

*Case* (i): $p < 1/2$.

(a) $\sum_{m=1}^{\infty} E(B_m) < \infty$ by virtue of (5.4).

(b) $\sum_{m=1}^{\infty} E(A_m)/m^{2-p+\varepsilon} = \sum_{m=1}^{\infty} E[A_m/m^{1-p+\varepsilon/2}]/m^{1+\varepsilon/2} < \infty$ by virtue of Lemma A.2 (in the Appendix).

(c) Clearly, $\sum_{m=1}^{\infty} c_{\varepsilon,p}/m^{2-p-\gamma_p\varepsilon} < \infty$.

*Case* (ii): $p = 1/2$. We need to divide both sides of (5.7) by $\log n$.

(a) $\sum_{m=1}^{n} E(B_m)/\log n \underset{n\to\infty}{\longrightarrow} p^2 c_p/2 = c_{1/2}/8$ by virtue of (5.4) since $E(B_m) = (p^2 c_p + \varepsilon_m)/2m$ where $\varepsilon_m \to 0$ as $m \to \infty$ and $\sum_{m=1}^{n} \frac{1}{m}/\log n \to 1$.

(b) $\sum_{m=1}^{\infty} E(A_m)/m^{2-p+\varepsilon} = \sum_{m=1}^{\infty} E[A_m/m^{1-p+\varepsilon/2}]/m^{1+\varepsilon/2} < \infty$ by virtue of Lemma A.2 (in the Appendix). Hence

$$\frac{\sum_{m=1}^{n} E(A_m)/m^{2-p+\varepsilon}}{\log n} \to 0 \qquad \text{as } n \to \infty.$$

(c) Clearly,

$$\frac{\sum_{m=1}^{n} c_{\varepsilon,p}/m^{2-p-\gamma_p\varepsilon}}{\log n} \underset{n\to\infty}{\longrightarrow} 0$$

(since the numerator is summable). Hence

$$\frac{E(A_n)}{n^{1/2}\log n} \underset{n\to\infty}{\longrightarrow} c_{1/2}/8.$$

*Case* (iii): $p > 1/2$. We need to divide both sides of (5.7) by $n^{2p-1}$.

(a) $\sum_{m=1}^{n} E(B_m)/n^{2p-1} \xrightarrow[n\to\infty]{} \frac{p^2}{2(2p-1)} c_p$ by virtue of (5.4) since $E(B_m) = (p^2 c_p/2 + \varepsilon_m)/m^{2-2p}$ where $\varepsilon_m \to 0$ as $m \to \infty$ and $\sum_{m=1}^{n} \frac{1}{m^{2-2p}}/n^{1-2p} \to \frac{1}{2p-1}$.

(b)

$$\frac{\sum_{m=1}^{n} E(A_m)/m^{2-p+\varepsilon}}{n^{2p-1}} = \frac{\sum_{m=1}^{n} E[(A_m)/m^{p+\varepsilon/2}]/m^{2-2p+\varepsilon/2}}{n^{2p-1}} \xrightarrow[n\to\infty]{} 0$$

by virtue of Lemma A.2 (in the Appendix).

(c) For small enough $\varepsilon$, $\frac{\sum_{m=1}^{n} c_{\varepsilon,p}/m^{2-p-\gamma p\varepsilon}}{n^{2p-1}} \xrightarrow[n\to\infty]{} 0$ (since the numerator is summable). Hence

$$\frac{E(A_n)}{n^p} \xrightarrow[n\to\infty]{} \frac{p^2}{2(2p-1)} c_p. \qquad \square$$

We now consider the almost sure convergence properties of the average rank of the items kept, suitably normalized. We shall need the following result, due to Robbins and Siegmund [7], quoted as follows:

PROPOSITION 5.1. *Let* $(\Omega, \mathcal{F}, P)$ *be a probability space and let* $\mathcal{F}_1 \subset \mathcal{F}_2 \subset \cdots$ *be a sequence of sub-$\sigma$-algebras of* $\mathcal{F}$. *For each* $n = 1, 2, \ldots,$ *let* $z_n$, $\beta_n$, $\xi_n$ *and* $\zeta_n$ *be nonnegative* $\mathcal{F}_n$-*measurable random variables such that*

$$E(z_{n+1}|\mathcal{F}_n) \leq z_n(1 + \beta_n) + \xi_n - \zeta_n.$$

*Then* $\lim_{n\to\infty} z_n$ *exists and is finite and* $\sum_{n=1}^{\infty} \zeta_n < \infty$ *a.s. on* $\{\sum_{n=1}^{\infty} \beta_n < \infty, \sum_{n=1}^{\infty} \xi_n < \infty\}$.

*The limiting behavior of* $A_n$, *the average rank of the retained observations, depends on* $p$. *Theorems* 5.2, 5.3 *and* 5.4 *show the results for* $0 < p < \frac{1}{2}, \frac{1}{2} < p \leq 1$ *and* $p = \frac{1}{2}$, *respectively.*

THEOREM 5.2. *If* $0 < p < \frac{1}{2}$, *then* $A_n/n^{1-p}$ *converges almost surely as* $n \to \infty$ *to a nondegenerate random variable.*

THEOREM 5.3. *If* $\frac{1}{2} < p \leq 1$, *then* $V_n \xrightarrow[n\to\infty]{a.s.} q_p$ *where* $V_n = Q_n/L_n^2 = A_n/L_n$ *and* $q_p = \frac{1}{2}p^2/(2p-1)$. *Furthermore,* $E(V_n) \to q_p$.

THEOREM 5.4. *If* $p = \frac{1}{2}$, *then* $V_n/\log n \xrightarrow[n\to\infty]{a.s.} 1/8$ *and* $E(V_n/\log n) \to 1/8$.

Proof of Theorem 5.2.   We show that the almost sure convergence of $A_n/n^{1-p}$ is the result of a direct application of Proposition 5.1 above. Regard (5.6). Note that $B_n$ of (5.3) can be written as $B_n = \frac{(p^2 L_n + \theta_n p)/2}{(n+1)^{2-p}}$ where $|\theta_n| \leq 1$. Apply Proposition 5.1 to (5.6) with $z_n = Y_n$, $\xi_n = B_n$, $\zeta_n = 0$, $\beta_n = \{e_n/[(L_n+1)(n+1)] + |O(n^{-2})|\}$, all nonnegative. Since $L_n/n^p$ converges a.s., $\sum_{n=1}^{\infty} \beta_n < \infty$ and $\sum_{m=1}^{\infty} \xi_m < \infty$ a.s. The nondegeneracy of the limit follows from the fact that the first observations have rank of order $n$ and their influence on $A_n/n^{1-p}$ does not vanish as $n \to \infty$.   $\square$

In order to obtain results for $p \geq 1/2$ we need two lemmas. These lemmas describe the extent to which the sums of the ranks of the items kept $Q_{n+1}$ increase after the $(n+1)$st item is observed. If the $(n+1)$st item is kept, it is easy to see that $Q_{n+1} = Q_n + L_n + 1$. The difficulty arises when the $(n+1)$st item is not kept. Even though the $(n+1)$st item is not kept, it might still be better than items that have been kept.

Let $X_i^*$ be the values of the items that are kept after having observed $n$ items for $i = 1, \ldots, L_n$ where the items are indexed from smallest (best) to largest (worst). The sum of the ranks of the kept items after $n+1$ items are observed, provided that the $(n+1)$st item is not kept, is denoted by $Q_n + A_{n+1}^*$ where $A_{n+1}^* = \sum_{i=j_n+1}^{L_n} I(X_{n+1} < X_i^*)$. The behavior of $A_{n+1}^*$ is captured as follows.

Lemma 5.1.

$$(5.8) \quad E(A_{n+1}^* | X_{n+1} > X_{j_n}^*, \mathcal{F}_n) = \frac{Q_n - (1/2)j_n(j_n+1) - j_n(L_n - j_n)}{n+1-j_n}$$

and

$$(5.9) \quad E(A_{n+1}^{*2} | X_{n+1} > X_{j_n}^*, \mathcal{F}_n) \leq \frac{Q_n(2L_n - 2j_n + 1)}{n+1-j_n},$$

where $j_n = \lceil pL_n \rceil$ as above.

Proof.

$$E(I(X_{n+1} < X_i^*) | X_{n+1} > X_{j_n}^*, \mathcal{F}_n) = \frac{R_i^{*^n} - j_n}{n+1-j_n}$$

for $i > j_n$ where $R_i^{*^n}$ is the rank of $X_i^*$ from among $X_1, \ldots X_n$

and

$$\sum_{i=j_n+1}^{L_n} (R_i^{*^n} - j_n) = Q_n - \tfrac{1}{2} j_n(j_n+1) - j_n(L_n - j_n),$$

hence (5.8) follows.

Since $I(X_{n+1} < X_i^*) = 1 \Rightarrow I(X_{n+1} < X_j^*) = 1$, for $j > i > j_n$,

$$E(A_{n+1}^{*2}|X_{n+1} > X_{j_n}^*, \mathcal{F}_n)$$

$$= E(A_{n+1}^*|X_{n+1} > X_{j_n}^*, \mathcal{F}_n)$$

$$+ 2E\left[\sum_{i=j_n+1}^{L_n-1}\sum_{k=i+1}^{L_n} I(X_{n+1} < X_i^*)|X_{n+1} > X_{j_n}^*, \mathcal{F}_n\right]$$

$$\leq E(A_{n+1}^*|X_{n+1} > X_{j_n}^*, \mathcal{F}_n)(2L_n - 2j_n + 1).$$

So (5.9) follows from (5.8). $\quad\square$

LEMMA 5.2.   Let $V_n = Q_n/L_n^2$. Then

$$(5.10)\quad E(V_{n+1}|\mathcal{F}_n) = V_n\left(1 - \frac{2p-1}{n+1} + O\left(\frac{1}{nL_n}\right)\right) + \frac{p^2/2}{n+1} + O\left(\frac{1}{nL_n}\right).$$

PROOF.

$$E(V_{n+1}|\mathcal{F}_n)$$

$$(5.11)\quad = \frac{j_n}{n+1}\frac{Q_n + L_n + 1}{(L_n+1)^2}$$

$$+ \frac{n+1-j_n}{(n+1)L_n^2}\left(Q_n + \frac{Q_n - (1/2)j_n(j_n+1) - j_n(L_n - j_n)}{n+1-j_n}\right).$$

But

$$\frac{j_n}{n+1}\frac{Q_n + L_n + 1}{(L_n+1)^2}$$

$$= \frac{j_n}{n+1}V_n\frac{L_n^2}{(L_n+1)^2} + \frac{j_n}{n+1}\frac{1}{L_n+1}$$

$$(5.12)\quad = \frac{j_n}{n+1}V_n\left(1 - \frac{1}{L_n+1}\right)^2 + \frac{j_n}{n+1}\frac{1}{L_n+1}$$

$$= \frac{j_n}{n+1}V_n + \frac{V_n}{n+1}\left(-2p + O\left(\frac{1}{L_n}\right)\right)$$

$$+ \frac{1}{n+1}\left(p + O\left(\frac{1}{L_n}\right)\right).$$

Also, the second term in (5.11) equals

$$(5.13)\quad \frac{n+1-j_n}{n+1}V_n + \frac{V_n}{n+1} - \frac{1}{n+1}\left(\frac{1}{2}p^2 + p(1-p) + O\left(\frac{1}{L_n}\right)\right)$$

$$= \frac{n+1-j_n}{n+1}V_n + \frac{V_n}{n+1} - \frac{1}{n+1}\left(p - \frac{1}{2}p^2 + O\left(\frac{1}{L_n}\right)\right).$$

Substituting (5.12) and (5.13) into (5.11) yields (5.10). □

LEMMA 5.3.

$$E(V_{n+1}^2|\mathcal{F}_n) = V_n^2 - \frac{2V_n^2}{n+1}\left(2p - 1 + O\left(\frac{1}{L_n}\right)\right)$$

(5.14)

$$+ \frac{V_n}{n+1}\left(p^2 + O\left(\frac{1}{L_n}\right)\right) + O\left(\frac{1}{nL_n}\right).$$

PROOF.

$$E(V_{n+1}^2|\mathcal{F}_n) = \frac{j_n}{n+1}\frac{(Q_n + L_n + 1)^2}{(L_n + 1)^4}$$

$$+ \frac{n+1-j_n}{(n+1)L_n^4}E[(Q_n + A_{n+1}^*)^2|X_{n+1} > X_{j_n}^*, \mathcal{F}_n].$$

But

$$\frac{(Q_n + L_n + 1)^2}{(L_n + 1)^4}$$

(5.15)

$$= \frac{Q_n^2}{L_n^4}\left(1 - \frac{1}{L_n + 1}\right)^4 + \frac{2Q_n}{(L_n + 1)^3} + O\left(\frac{1}{L_n^2}\right)$$

$$= V_n^2\left(1 - \frac{4}{L_n + 1} + O\left(\frac{1}{L_n^2}\right)\right)$$

$$+ \frac{2Q_n}{(L_n + 1)^3} + O\left(\frac{1}{L_n^2}\right).$$

Hence,

$$\frac{j_n}{n+1}\frac{(Q_n + L_n + 1)^2}{(L_n + 1)^4}$$

(5.16)

$$= \frac{j_n}{n+1}V_n^2 - \frac{V_n^2}{n+1}\left(4p + O\left(\frac{1}{L_n}\right)\right)$$

$$+ \frac{V_n}{n+1}\left(2p + O\left(\frac{1}{L_n}\right)\right) + O\left(\frac{1}{nL_n}\right).$$

Similarly, by Lemma 5.1,

$$\frac{n+1-j_n}{(n+1)L_n^4}E[(Q_n + A_{n+1}^*)^2|X_{n+1} > X_{j_n}^*, \mathcal{F}_n]$$

$$= \frac{n+1-j_n}{n+1}V_n^2 + \frac{2Q_n}{(n+1)L_n^4}\left(Q_n - \frac{1}{2}j_n(j_n + 1) - j_n(L_n - j_n)\right)$$

(5.17)

$$+ \frac{n+1-j_n}{(n+1)L_n^4} E(A_{n+1}^{*2}|X_{n+1} > X_{j_n}^*, \mathcal{F}_n)$$

$$= \frac{n+1-j_n}{(n+1)} V_n^2 + \frac{2V_n^2}{n+1} - \frac{2V_n}{n+1}\left(\frac{1}{2}p^2 + p(1-p) + O\left(\frac{1}{L_n}\right)\right)$$

since by (5.9) the term involving $E(A_{n+1}^{*2}|X_{n+1} > X_{j_n}^*, \mathcal{F}_n)$ is bounded by

$$\frac{Q_n(2L_n - 2j_n + 1)}{(n+1)L_n^4} = \frac{V_n}{n+1}O\left(\frac{1}{L_n}\right).$$

Combining (5.16) and (5.17) yields (5.14).  $\square$

PROOF OF THEOREM 5.3.   The claim is trivial for $p = 1$, for then $Q_n = \frac{1}{2}L_n(L_n + 1)$. Suppose $\frac{1}{2} < p < 1$. Let $V_n^* = (V_n - q_p)^2$. We shall show that $\check{V}_n^* \to 0$ a.s.:

$$E(V_{n+1}^*|\mathcal{F}_n) = E(V_{n+1}^2|\mathcal{F}_n) - 2q_p E(V_{n+1}|\mathcal{F}_n) + q_p^2.$$

Lemmas 5.2 and 5.3 imply that

$$\begin{aligned}
(5.18) \quad E(V_{n+1}^*|\mathcal{F}_n) = {} & V_n^* - \frac{2V_n^2}{n+1}\left(2p - 1 + O_1\left(\frac{1}{L_n}\right)\right) \\
& + \frac{V_n}{n+1}\left(p^2 + 2q_p(2p-1) + O_2\left(\frac{1}{L_n}\right)\right) \\
& - \frac{q_p p^2}{n+1}\left(1 + O_3\left(\frac{1}{L_n}\right)\right).
\end{aligned}$$

Notice that $|V_n - q_p| = \sqrt{V_n^*}$ so that $V_n \le \sqrt{V_n^*} + q_p$. Hence $0 \le V_n \le V_n^* + 1 + q_p$.

$V_n^2 = V_n^* + 2q_p V_n - q_p^2$. Recalling that $q_p = \frac{1}{2}p^2/(2p-1)$, we obtain [after algebra, applying (5.18)]

$$\begin{aligned}
(5.19) \quad E(V_{n+1}^*|\mathcal{F}_n) = {} & \left(1 - \frac{2(2p-1)}{n+1}\right)V_n^* + \frac{2(2p-1)}{n+1}O_1\left(\frac{1}{L_n}\right) \\
& + \frac{V_n}{n+1}O_2\left(\frac{1}{L_n}\right) + \frac{q_p p^2}{n+1}O_3\left(\frac{1}{L_n}\right) \\
\le {} & \left(1 - \frac{2(2p-1) + O_2(1/L_n)}{n+1}\right)V_n^* + O\left(\frac{1}{nL_n}\right).
\end{aligned}$$

Since $L_n/n^p \xrightarrow{a.s.}$ to a finite random variable, there exists (a random) $n_0$ such that $|O_2(\frac{1}{L_{n_0}})| < (2p - 1)$ for all $n \ge n_0$, so that for all $n \ge n_0$,

$$(5.20) \quad E(V_{n+1}^*|\mathcal{F}_n) \le V_n^* + \left|O\left(\frac{1}{nL_n}\right)\right|.$$

But since $L_n/n^p$ converges a.s., it follows that $|O(\frac{1}{nL_n})|$ is a.s. summable. Therefore, it follows by Proposition 5.1 that $V_n^*$ converges almost surely.

We claim that $V_n^* \xrightarrow{a.s.} 0$. Suppose this is false and there exists a set $A$ of positive measure such that $\lim V_n^* \geq a$ a.s. on $A$, where $a > 0$. Since $\lim V_n^*$ exists, (5.19) can be rewritten as

$$(5.21) \qquad E(V_{n+1}^* | \mathcal{F}_n) \leq \left(1 - \frac{2(2p-1)}{n+1}\right) V_n^* + O\left(\frac{1}{nL_n}\right).$$

For given $0 < \varepsilon < a$, there exists (a random) $n_1 \geq n_0$ such that, on $A$, $V_n^* > a - \varepsilon$ for all $n \geq n_1$. Clearly $P(A | \mathcal{F}_{n_1}) > 0$. From (5.21) it follows that

$$E(V_n^* | \mathcal{F}_{n_1}) = \sum_{i=n_1+1}^{n} E(V_i^* - V_{i-1}^* | \mathcal{F}_{n_1}) + V_{n_1}^*$$

$$\leq -2(2p-1)(a-\varepsilon) P(A | \mathcal{F}_{n_1}) \sum_{i=n_1+1}^{n} \frac{1}{i+1} + V_{n_1}^* + c$$

$$\xrightarrow[n \to \infty]{} -\infty$$

where $c = \sum_{i=n_1+1}^{\infty} |O(\frac{1}{iL_i})| < \infty$.

This is clearly a contradiction since $V_n^*$ is nonnegative. Therefore, $V_n^* \to 0$ a.s., that is, $V_n \to q_p$ a.s.

By virtue of Lemma A.2 in the Appendix, there exists $\eta > 0$ such that $E(\frac{1}{nL_n}) < \frac{1}{n^{1+\eta}}$. Therefore it follows from (5.20) that $\{E(V_n)^*\}$ is bounded. $V_n \to q_p$ a.s. now implies $EV_n \xrightarrow[n \to \infty]{} q_p$. $\square$

SKETCH OF PROOF OF THEOREM 5.4.    Let $V_n^* = (\frac{V_n}{\log n} - \frac{1}{8})^2$.

Ignoring $O(\cdot)$ terms, it follows from (5.10) that

$$E\left(\frac{V_{n+1}}{\log(n+1)} \Big| \mathcal{F}_n\right) \approx \frac{V_n}{\log n} \frac{\log n}{\log(n+1)} + \frac{1/8}{n \log(n+1)}$$

$$(5.22) \qquad \approx \frac{V_n}{\log n}\left(1 - \frac{1}{n \log(n+1)}\right) + \frac{1/8}{n \log(n+1)}.$$

The last expression comes from $\log n / \log(n+1) = 1 - \frac{1}{n \log(n+1)} + O(\frac{1}{n^2 \log n})$.

Similarly, from (5.14),

$$E\left(\frac{V_{n+1}^2}{\log^2(n+1)} \Big| \mathcal{F}_n\right)$$

$$(5.23) \qquad \approx \frac{V_n^2}{\log^2 n}\left(1 - \frac{1}{n \log(n+1)}\right)^2$$

$$+ \frac{1}{4} \frac{V_n}{\log n} \cdot \frac{1}{n \log(n+1)} \left( 1 - \frac{1}{n \log(n+1)} \right)$$

$$\approx \frac{V_n^2}{\log^2 n} - \frac{2 V_n^2}{\log^2 n} \frac{1}{n \log(n+1)} + \frac{1}{4} \frac{V_n}{\log n} \frac{1}{n \log(n+1)}.$$

But

$$(5.24) \quad E(V_{n+1}^* | \mathcal{F}_n) = E\left[ \frac{V_{n+1}^2}{\log^2(n+1)} \Big| \mathcal{F}_n \right] - \frac{1}{4} E\left[ \frac{V_{n+1}}{\log(n+1)} \Big| \mathcal{F}_n \right] + \frac{1}{64}.$$

Substituting (5.22) and (5.23) into (5.24) yields

$$E(V_{n+1}^* | \mathcal{F}_n) \approx V_n^* - \frac{2}{n \log(n+1)} \left( \frac{V_n}{\log n} - \frac{1}{8} \right)^2$$

$$= V_n^* - \frac{2}{n \log(n+1)} V_n^*.$$

The remainder of the proof follows in a fashion similar to the proof of Theorem 5.3. □

**6. Discussion.** An extensive simulation study has been performed with varying horizons up to 10,000 and also 10,000 replications. The results of this simulation will be published elsewhere. The estimated constants for the limit of $E(L_n/n^p)$ (see Theorem 4.1) and $E(A_n/a_n(p))$ (see Theorem 5.1) are presented in Table 1.

The standard errors for estimation of the constants are less than 0.002.

REMARK 6.1. Some of the results of the preceding section carry over to the inverse problem of fixing the number of items kept and considering the number of observations required until this goal is achieved. Suppose that a $p$-percentile rule is applied. Let $Z_i$ (for $i \geq 1$) be the number of observations made from the instant that the size of the set of retained observations became $i - 1$ until its size became $i$. Also, let $N_n = \sum_{i=1}^n Z_i$ be the number of observations made until $n$ have been retained. The results stated in Theorems 4.3 and 5.2, 5.3 and 5.4 carry over directly to $N_n$ and $Q_{N_n}: L_{N_n}/N_n^p = n/N_n^p$ converges a.s. as $n \to \infty$ to a finite, nondegenerate random variable. For $p < \frac{1}{2}$, $p = \frac{1}{2}$, $p > \frac{1}{2}$ the quantities

TABLE 1
$a_n(p) = n^{1-p}$ if $p < 0.5$, $a_n(\frac{1}{2}) = n^{1/2} \log n$, $a_n(p) = n^p$ if $p > 0$

| p | 0.1 | 0.2 | 0.3 | 0.4 | 0.5 | 0.6 | 0.7 | 0.8 | 0.9 | 1.0 |
|---|---|---|---|---|---|---|---|---|---|---|
| $E(L_n/n^p)^a$ | 4.178 | 2.674 | 2.111 | 1.653 | 1.181 | 1.198 | 1.045 | 0.841 | 0.693 | 0.500 |
| $E(A_n/a_n)^a$ | 0.238 | 0.401 | 0.578 | 0.967 | 0.214 | 0.978 | 0.634 | 0.449 | 0.351 | 0.250 |

[a]Estimate of the limit using $n = 10{,}000$ and 10,000 replications.

($\alpha$)  $\frac{Q_{N_n}}{n}/N_n^{1-p}, \frac{Q_{N_n}}{n}/(N_n^{1/2}\log N_n), \frac{Q_{N_n}}{n}/N_n^p$,

respectively, converge a.s. as $n \to \infty$ to nondegenerate random variables, and for $p > \frac{1}{2}$,

($\beta$)  $Q_{N_n}/n^2 \underset{n\to\infty}{\longrightarrow} p^2/[2(2p-1)]$  a.s.

For $0 < p \le 1$ and $n \ge 2$ the expectation of $N_n$ is $\infty$ [since $E(Z_2) = \infty$]. For $p < \frac{1}{2}, p = \frac{1}{2}, p > \frac{1}{2}$, respectively, the corresponding expectation of each of the three expressions in ($\alpha$) converges to finite positive constants as $n \to \infty$, and for $p > \frac{1}{2}$ the expectation of the term on the left-hand side of ($\beta$) converges to the value in the right side, as $n \to \infty$.

REMARK 6.2.   Define $R_i^n = 1 + \sum_{j=1}^n I(X_j < X_i)$. (This does not constitute a change when the $X$'s are i.i.d. continuous.) Then Theorem 2.1 and its proof hold verbatim, even if there are ties among the $X$'s (which may be the case if the $X$'s are discrete random variables), and even if the $X$'s are not random.

REMARK 6.3.   Other rules can be evaluated in a manner similar to the $p$-percentile rules. For example, for $k$-record rules using Lemma 2.1, it can be shown that $L_n/\log n$ converges almost surely to $k$ as $n \to \infty$ and $E(L_n)/\log n$ also converges to $k$. It can be shown that $Q_n/(n+1)$ is a (nonnegative) submartingale that converges a.s. as $n \to \infty$ to a nondegenerate random variable, and $E(Q_n)/(n+1) \underset{n\to\infty}{\longrightarrow} k$. Thus $A_n \log n/n$ converges a.s. to a nondegenerate random variable.

REMARK 6.4.   The complexity of the sorting problem is of order $n \log n$. Sometimes, one is interested in retaining (in sorted form) only some of the best observations rather than the whole set. In this case, a $p$-percentile rule with $0 < p < 1$ obtains a sorted set of best observations, and the complexity is of order $n$. To see this, note that initially each observation has to be compared only to the $p$-percentile of the retained set—amounting to $n$ operations—and each retained observation must be compared to (roughly) $100p\%$ of the retained observations, amounting (at most) to another $2\sum_{j=1}^{L_n} \log(j+1) = O(L_n \log L_n) = O_p(n^p \log n)$ operations (since the retained set can be stored in sorted condition).

## APPENDIX

LEMMA A.1.   *For any $p$-percentile rule with $0 < p \le 1$ and any $0 < \varepsilon < 1$, there exists a constant $0 < c_{\varepsilon,p} < \infty$ such that*

$$P(L_n < k) \le c_{\varepsilon,p} k^{r_0(1-\varepsilon)}/n^{1-\varepsilon}$$

*where $r_0 = \lceil 1/p \rceil$.*

PROOF.    After $m-1$ observations have been retained by the $p$-percentile rule, let $Z_m$ denote the number of additional observations until the next retention. Note that $Z_1 = 1$. Let $N_n = \sum_{i=1}^{n} Z_i$ be the number of observations made until $n$ items have been retained. Thus,

$$(A.1) \quad P(L_n < k) = P(N_k > n) = P(N_k^{1-\varepsilon} > n^{1-\varepsilon}) \leq E(N_k^{1-\varepsilon})/n^{1-\varepsilon}.$$

Without loss of generality, assume that the $X_i$ have a $U[0,1]$ distribution. Let $X_i^n$ denote the observation with rank $i$ among $X_1, X_2, \ldots, X_n$. Note that conditional on $X_1, \ldots, X_{N_m}$, the distribution of $Z_{m+1}$ is

$$\text{Geometric } p \text{ with } p = X_{j_{N_m}}^{N_m} = X_{\lceil pm \rceil}^{N_m}.$$

Also note that conditional on $N_m$, the distribution of $X_{\lceil pm \rceil}^{N_m}$ is Beta $(\lceil pm \rceil, N_m + 1 - \lceil pm \rceil)$.

Therefore, for $0 \leq \varepsilon < 1$

$$
\begin{aligned}
E(Z_{m+1}^{1-\varepsilon}|N_m) &= E[E(Z_{m+1}^{1-\varepsilon}|N_m, X_{\lceil pm \rceil}^{N_m})|N_m] \\
&\leq E([E(Z_{m+1}|N_m, X_{\lceil pm \rceil}^{N_m})]^{1-\varepsilon}|N_m) \\
&= E\left[\left(\frac{1}{X_{\lceil pm \rceil}^{N_m}}\right)^{1-\varepsilon}\Big|N_m\right] \\
(A.2) \qquad &= \frac{N_m!}{(\lceil pm \rceil - 1)!(N_m - \lceil pm \rceil)!} \\
&\qquad \times \frac{\Gamma(\lceil pm \rceil + \varepsilon - 1)\Gamma(N_m + 1 - \lceil pm \rceil)}{\Gamma(N_m + \varepsilon)} \\
&= \frac{N_m!}{(\lceil pm \rceil - 1)!} \cdot \frac{\Gamma(\lceil pm \rceil + \varepsilon - 1)}{\Gamma(N_m + \varepsilon)} \\
&\leq \frac{N_m}{\lceil pm \rceil + \varepsilon - 1}.
\end{aligned}
$$

The last inequality in (A.2) follows since $f(A) = \Gamma(A + \varepsilon)/(\Gamma(A))$ is increasing in $A$ for any integer $A$. This is easily seen since $f(A+1)/f(A) = (A+\varepsilon)/A \geq 1$ for all $\varepsilon \geq 0$.

For $\varepsilon = 0$ and $\lceil pm \rceil > 1$ obtain $E(Z_{m+1}|N_m) \leq N_m/(\lceil pm \rceil - 1)$, so that

$$E(N_{m+1}|N_m) \leq \frac{\lceil pm \rceil}{\lceil pm \rceil - 1} N_m.$$

Hence

$$E(N_{m+1}^{1-\varepsilon}|N_m) \leq [E(N_{m+1}|N_m)]^{1-\varepsilon} \leq \left(\frac{\lceil pm \rceil}{\lceil pm \rceil - 1}\right)^{1-\varepsilon} N_m^{1-\varepsilon}.$$

Letting $m_p$ be the smallest $m$ such that $\lceil pm \rceil > 1$, it follows that

$$(A.3) \qquad E(N_{m+1}^{1-\varepsilon} | N_{m_p}) \le N_{m_p}^{1-\varepsilon} \prod_{i=m_p}^{m} \left( \frac{\lceil pi \rceil}{\lceil pi \rceil - 1} \right)^{1-\varepsilon}.$$

We first show that $E(N_{m_p}^{1-\varepsilon})$ is finite. By virtue of (A.2), since $Z_{m+1} \ge 1$,

$$\begin{aligned}
E(N_{m+1}^{1-\varepsilon} | N_m) &= E((N_m + Z_{m+1})^{1-\varepsilon} | N_m) \\
&\le E(N_m^{1-\varepsilon} | N_m) + E(Z_{m+1}^{1-\varepsilon} | N_m) \\
&\le N_m^{1-\varepsilon} + E(Z_{m+1}^{(1-\varepsilon/2)^2} | N_m) \\
&\le N_m^{1-\varepsilon} + [E(Z_{m+1}^{1-\varepsilon/2} | N_m)]^{1-\varepsilon/2} \\
&\le N_m^{1-\varepsilon} + \frac{N_m^{1-\varepsilon/2}}{(\lceil pm \rceil + \varepsilon/2 - 1)^{1-\varepsilon/2}}.
\end{aligned}$$

This recursive relation can be applied repeatedly. Since $E(N_1^{1-\varepsilon}) = 1$, it follows that $E(N_{m_p}^{1-\varepsilon})$ is finite for all $0 < \varepsilon < 1$.

Finally, note that $\frac{\lceil pi \rceil}{\lceil pi \rceil - 1} > 1$ and can appear (in the product $\prod_{i=m_p}^{m} \frac{\lceil pi \rceil}{\lceil pi \rceil - 1}$) at most $r_0 = \lceil 1/p \rceil$ times. Hence,

$$\prod_{i=m_p}^{m} \left( \frac{\lceil pi \rceil}{\lceil pi \rceil - 1} \right)^{1-\varepsilon} \le \left( \frac{\lceil pm \rceil}{\lceil pm_p \rceil - 1} \right)^{r_0(1-\varepsilon)} = \lceil pm \rceil^{r_0(1-\varepsilon)}.$$

For $m = k$, the conclusion follows from (A.1) and (A.3). $\quad\square$

LEMMA A.2. *Let* $0 < \varepsilon < 1$ *and consider a $p$-percentile rule with* $0 < p \le 1$. *Let* $p^* = \max(p, 1 - p)$. *Then*

$$\lim_{n \to \infty} E\left( \frac{A_n}{n^{p^* + \varepsilon}} \right) = 0.$$

PROOF. The statement is trivially true for $p = 1$ as $E(A_n/n) \to 1/4$ as $n \to \infty$. Thus consider $0 < p < 1$. Let $1 \le k_\varepsilon < \infty$ be an integer such that

$$\frac{1 + L_n - j_n}{1 + L_n} < (1-p)(1+\varepsilon) \le p^*(1+\varepsilon)$$

whenever $L_n \ge k_\varepsilon$, where $j_n = \lceil pL_n \rceil$.

From Lemma A.1 we obtain

$$P(L_n < k_\varepsilon) \le c^*_{\varepsilon, p, k_\varepsilon} / n^{1-\varepsilon}.$$

Using Lemma 2.1(iii),

$$E(A_{n+1}|\mathcal{F}_n) = A_n\left(1 + \frac{1 + L_n - j_n}{(n+1)(1+L_n)}\right) + \frac{j_n(j_n-1)/2}{(n+1)L_n}$$

$$\text{(A.4)} \qquad \leq A_n + A_n\frac{1+L_n-j_n}{(n+1)(1+L_n)}[I\{L_n \geq k_\varepsilon\} + I\{L_n < k_\varepsilon\}]$$

$$+ \frac{p^2 L_n + p}{2(n+1)}.$$

Since $A_n \leq n$, the right-hand side of (A.4) is less than or equal to

$$A_n\left[1 + \frac{p^*(1+\varepsilon)}{n}\right] + \frac{n}{n+1}I\{L_n < k_\varepsilon\} + \frac{p^2 L_n + p}{2(n+1)}.$$

Hence, with $\varepsilon < \{p \wedge c_p\}$ and large enough $n$,

$$E(A_{n+1}) \leq E(A_n)\left[1 + \frac{p^*(1+\varepsilon)}{n}\right] + \frac{c^*_{\varepsilon,p,k_\varepsilon}}{n^{1-\varepsilon}} + \frac{p^2(\varepsilon + c_p)}{2n^{1-p}}$$

$$\leq E(A_n)\left[1 + \frac{p^*(1+\varepsilon)}{n}\right] + \frac{c_p}{n^{1-p}}.$$

Let $\gamma_1 = 1$ and define $\gamma_{n+1} = \gamma_n/[1 + \frac{p^*(1+\varepsilon)}{n}]$. Thus, $\{\gamma_n\}$ is a decreasing sequence, and there exists $0 < \gamma_\infty < \infty$ such that

$$\lim_{n\to\infty} n^{p^*(1+\varepsilon)}\gamma_n = \gamma_\infty.$$

Note that for large enough $n$,

$$\gamma_{n+1}E(A_{n+1}) \leq \gamma_n E(A_n) + \frac{2c_p\gamma_\infty}{n^{1-p+p^*(1+\varepsilon)}}.$$

Because $1 - p + p^* \geq 1$, it follows that

$$\limsup_{n\to\infty} \gamma_n E(A_n) < \infty, \qquad \text{that is, } \limsup_{n\to\infty} E\left(\frac{A_n}{n^{p^*(1+\varepsilon)}}\right) < \infty.$$

Hence

$$\lim_{n\to\infty} E\left(\frac{A_n}{n^{p^*+\varepsilon}}\right) = 0. \qquad\qquad \square$$

**Acknowledgment.** The suggestions from a referee helped greatly to improve the exposition.

A. M. KRIEGER
DEPARTMENT OF STATISTICS
THE WHARTON SCHOOL
UNIVERSITY OF PENNSYLVANIA
PHILADELPHIA, PENNSYLVANIA 19104
USA
E-MAIL: krieger@wharton.upenn.edu

M. POLLAK
E. SAMUEL-CAHN
DEPARTMENT OF STATISTICS
HEBREW UNIVERSITY
MOUNT SCOPUS, JERUSALEM 91905
ISRAEL
E-MAIL: msmp@mscc.huji.ac.il
        scahn@mscc.huji.ac.il